\documentclass[11pt]{article}

\usepackage{amsmath, amssymb, amscd}

\def\eqref#1{(\ref{#1})}
\newcommand{\goth}{\mathfrak}

\newcommand{\arrow}{{\:\longrightarrow\:}}
\newcommand{\Z}{{\Bbb Z}}
\newcommand{\C}{{\Bbb C}}
\newcommand{\R}{{\Bbb R}}
\newcommand{\Q}{{\Bbb Q}}

\newcommand{\6}{\partial}
\def\1{\sqrt{-1}\:}
\newcommand{\restrict}[1]{{\left|_{{\phantom{|}\!\!}_{#1}}\right.}}
\newcommand{\cntrct}                
{\hspace{2pt}\raisebox{1pt}{\text{$\lrcorner$}}\hspace{2pt}}

\def\Bbb#1{\mathbb #1}

\renewcommand{\bar}{\overline}
\renewcommand{\phi}{\varphi}
\renewcommand{\epsilon}{\varepsilon}

\newcommand{\End}{\operatorname{End}}

\newcommand{\Id}{\operatorname{Id}}

\newcommand{\Vol}{\operatorname{Vol}}
\newcommand{\Hom}{\operatorname{Hom}}

\newcommand{\Alt}{\operatorname{Alt}}

\newcommand{\Tw}{\operatorname{Tw}}
\newcommand{\Tr}{\operatorname{Tr}}

\renewcommand{\Re}{\operatorname{Re}}
\renewcommand{\Im}{\operatorname{Im}}

\newcommand{\comment}[1]{{}}

\def\blacksquare{\hbox{\vrule width 4pt height 4pt depth 0pt}}
\def\endproof{\blacksquare}

\makeatletter

\newcommand{\ps@verbit}{%
  \renewcommand{\@oddhead}{%
          \scriptsize
          {Almost complex 6-manifolds}
          \hfil\tiny {M. Verbitsky, July 7, 2005}}
  \renewcommand{\@evenhead}{\@oddhead}
  \renewcommand{\@oddfoot}{\hfil\thepage\hfil}
  \renewcommand{\@evenfoot}{\@oddfoot}}
 
\newcounter{Mycounter}[section]
\newcounter{lemma}[section]
\renewcommand{\thelemma}{\noindent{Lemma \thesection.\arabic{lemma}}}
\newcommand{\lemma}{%
     \setcounter{lemma}{\value{Mycounter}}
     \refstepcounter{lemma}
     \stepcounter{Mycounter}
     {\bf \thelemma:\ }}

\newcounter{claim}[section]
\renewcommand{\theclaim}{\noindent{Claim \thesection.\arabic{claim}}}
\newcommand{\claim}{%
     \setcounter{claim}{\value{Mycounter}}
     \refstepcounter{claim}
     \stepcounter{Mycounter}
     {\bf \theclaim:\ }}

\newcounter{sublemma}[section]
\newcounter{corollary}[section]
\renewcommand{\thecorollary}{\noindent{Corollary \thesection.\arabic{corollary}}}
\newcommand{\corollary}{%
     \setcounter{corollary}{\value{Mycounter}}
     \refstepcounter{corollary}
     \stepcounter{Mycounter}
     {\bf \thecorollary:\ }}

\newcounter{theorem}[section]
\renewcommand{\thetheorem}{\noindent{Theorem \thesection.\arabic{theorem}}}
\newcommand{\theorem}{%
     \setcounter{theorem}{\value{Mycounter}}
     \refstepcounter{theorem}
     \stepcounter{Mycounter}
     {\bf \thetheorem:\ }}

\newcounter{conjecture}[section]
\newcounter{proposition}[section]
\renewcommand{\theproposition}
       {\noindent{Proposition \thesection.\arabic{proposition}}}
\newcommand{\proposition}{%
     \setcounter{proposition}{\value{Mycounter}}
     \refstepcounter{proposition}
     \stepcounter{Mycounter}
     {\bf \theproposition:\ }}

\newcounter{definition}[section]
\renewcommand{\thedefinition}
       {\noindent{Definition~\thesection.\arabic{definition}}}
\newcommand{\definition}{%
     \setcounter{definition}{\value{Mycounter}}
     \refstepcounter{definition}
     \stepcounter{Mycounter}
     {\bf \thedefinition:\ }}

\newcounter{example}[section]
\newcounter{remark}[section]
\renewcommand{\theremark}{\noindent{Remark \thesection.\arabic{remark}}}
\newcommand{\remark}{%
     \setcounter{remark}{\value{Mycounter}}
     \refstepcounter{remark}
     \stepcounter{Mycounter}
     {\bf \theremark:\ }}

\newcounter{problem}[section]
\newcounter{question}[section]
\@addtoreset{equation}{section}
\@addtoreset{footnote}{section}
\makeatother

\begin{document}

\begin{center}
{\LARGE\bf An intrinsic volume 
functional \\[2mm] on almost complex 6-manifolds\\[3mm] and
nearly K\"ahler geometry}
\\[4mm]
Misha Verbitsky\footnote{Misha Verbitsky is 
an EPSRC advanced fellow 
supported by CRDF grant RM1-2354-MO02 
and EPSRC grant  GR/R77773/01}
\\[4mm]

{\tt verbit@maths.gla.ac.uk, \ \  verbit@mccme.ru}
\end{center}

{\small 
\hspace{0.15\linewidth}
\begin{minipage}[t]{0.7\linewidth}
{\bf Abstract} \\
Let $(M,I)$ be an almost complex 6-manifold. 
The obstruction to integrability of the almost
complex structure (the so-called Nijenhuis tensor)
$N:\; \Lambda^{0,1}(M) \arrow \Lambda^{2,0}(M)$
maps a 3-dimensional bundle to a 3-dimensional one. 
We say that Nijenhuis tensor is {\bf non-degenerate}
if it is an isomorphism. An almost complex manifold 
$(M,I)$ is called {\bf nearly K\"ahler} if 
it admits a Hermitian form $\omega$ such that
$\nabla(\omega)$ is totally antisymmetric,
$\nabla$ being the Levi-Civita connection.
We show that a nearly K\"ahler metric on
a given almost complex 6-manifold with
non-degenerate Nijenhuis tensor is unique
(up to a constant). We interpret the nearly K\"ahler
property in terms of $G_2$-geometry and
in terms of connections with totally
antisymmetric torsion, obtaining
a number of equivalent definitions.

\ \ \ We construct a natural diffeomorphism-invariant
functional $I \arrow \int_M \Vol_I$ on the 
space of almost complex structures on $M$,
similar to the Hitchin functional, and compute
its extrema in the following important case. 
Consider an almost complex structure 
$I$ with non-degenerate Nijenhuis tensor,
admitting a Hermitian connection with totally
antisymmetric torsion. We show that
the Hitchin-like functional $I \arrow \int_M \Vol_I$
has an extremum in $I$ if and only if $(M,I)$ is nearly 
K\"ahler. 

\end{minipage}
}

{
\small
\tableofcontents
}


\setcounter{section}{-1}
\section{Introduction}
\label{_Intro_real_Section_}


\subsection{Almost complex manifolds with non-degenerate Nijenhuis tensor}

In geometry, two kinds of plane distributions often arise.
There are integrable ones: complex structures, foliations,
CR-structures. On the other hand, there are
``maximally non-integrable'' distributions, such as
the contact structures, where the obstruction
to integrability is nowhere degenerate. Looking 
at almost complex structures in dimension 3, 
one finds that the obstruction to the
integrability, so-caled Nijenhuis tensor
\[
N:\; \Lambda^2 T^{1,0}(M) \arrow T^{0,1}(M)
\]
maps a 3-dimensional bundle to a 3-dimensional bundle.
It is only natural to study the class of complex 3-manifolds
such that $N$ is nowhere degenerate. 

Given such a manifold $M$, it is possible to 
construct a nowhere degenerate, positive volume form  
$\det N^* \otimes \overline{\det N^*}$ on $M$ 
(for details, see \eqref{_volume_via_N_Equation_}).

We study extrema of this volume form, showing that these
extrema correspond to an interesting geometric structure
(\ref{_extre_for_M_I_Theorem_}).

In Hermitian geometry, one often encounters 
a special kind of almost complex Hermitian manifolds,
called strictly nearly K\"ahler (NK-) manifolds, or Gray
manifolds, after  Alfred Gray (\ref{_NK_tens_Definition_}). 
These manifolds can be characterized in terms of the $G_2$-structure
on their Riemannian cone, or in terms of a special set of equations
reminiscent of Calabi-Yau equations (Subsection 
\ref{_Nearly_Kah_G_2_Subsection_}).

We prove that a strictly nearly K\"ahler 3-manifold 
is uniquely determined by its almost complex structure 
(\ref{_NK_unique_from_I_Corollary_}). Moreover,
such manifolds are extrema of the volume
functional associated with the Nijenhuis tensor
(\ref{_extre_for_M_I_Theorem_}). This reminds
of a construction of Hitchin's functional
on the space of all $SL(3, \C)$-structures
on a manifold, having extrema on 
Calabi-Yau manifolds (\cite{_Hitchin:3-forms_}).

\subsection{Contents}

This paper has the following structure.

\hfill

$\bullet$ In Section \ref{_Intro_Section_},
we introduce the class of 3-manifolds with nowhere
degenerate Nijenhuis tensor, and describe the basic
structures associated with these manifolds. 
We give a sketch of a proof of existence of
a Hermitian connection with totally antisymmetric
torsion, due to Friedrich and Ivanov, and show that 
such a Hermitian metric is uniquely determined
by the almost complex structure, if the 
Nijenhuis tensor is nowhere degenerate. 

$\bullet$ In Section \ref{_NK_Intro_Section_},
we introduce the nearly K\"ahler manifolds, giving several
versions of their definition and listing examples.

$\bullet$ In Section \ref{_conne_antisymme_Nije_Section_}
we apply the results about connections with totally
antisymmetric torsion to nearly K\"ahler geometry, showing
that an almost complex structure determines the Hermitian
structure on such a manifold uniquely, up to a constant
multiplier. 

$\bullet$ In Section \ref{_NK_and_G2_antisy_Section_},
we give several additional versions of a definition
of a nearly K\"ahler manifold, obtaining an explicit
description of a Nijenhuis tensor in terms of
an orthonormal frame. We also interpret the 
nearly K\"ahler structure on a manifold
in terms of $G_2$-geometry of its 
Riemannian cone. This is used to show that
an NK-structure on a manifold $M$ is uniquely determined by
its metric, unless $M$ is locally isometric to a 
6-sphere (\ref{_NK_determined_by_metric_Proposition_}).

$\bullet$ In Section \ref{_varia_NK_Section_},
we study infinitesimal variations of an
almost complex structure. We prove that
NK-manifolds are extrema of an intrinsic
 volume functional described earlier.
A partial converse result is also obtained.
Given an almost complex manifold $M$ 
with nowhere degenerate Nijenhuis tensor,
admitting a Hermitian connection with
totally antisymmetric torsion, $M$
is an extremum of the intrinsic
volume functional if and only if
$M$ is nearly K\"ahler.


\section{Almost complex manifolds
with non-degenerate Nijenhuis tensor}
\label{_Intro_Section_}


\subsection{Nijenhuis tensor on 6-manifolds}
\label{_Nije_non-dege_Subsection_}

Let $(M,I)$ be an almost complex manifold.
The Nijenhuis tensor maps two $(1,0)$-vector fields
to the $(0,1)$-part of their commutator. This map is
$C^\infty$-linear, and vanishes, as the Newlander-Nirenberg
theorem implies, precisely when $I$ is integrable.
We write the Nijenhuis tensor as
\[
N:\; \Lambda^2 T^{1,0}(M)\arrow T^{0,1}(M).
\]
The dual map
\begin{equation}\label{_Nije_dual_Equation_}
N^*:\; \Lambda ^{0,1}(M) \arrow \Lambda^{2,0}(M)
\end{equation}
is also called the Nijenhuis tensor.
Cartan's formula implies that $N^*$ acts on $\Lambda^1(M)$
as the $(2, -1)$-part of the de Rham differential.

When one studies the distributions, one is usually interested
in integrable ones (such as $T^{1,0}(M)\subset TM\otimes
\C$  for complex or CR-manifolds) or ones where the obstruction
to integrability is nowhere degenerate (such as a contact 
distribution). 

For the Nijenhuis tensor in complex dimension $>3$,
non-degeneracy does not make much sense, because the space
$\Hom(\Lambda ^{0,1}(M), \Lambda^{2,0}(M))$ becomes quite
complicated. However, for $n=3$, both sides of
\eqref{_Nije_dual_Equation_} are 3-dimensional,
and we can define the non-degeneracy as follows.

\hfill

\definition
Let $(M,I)$ be an almost complex manifold of real
dimension 6, and $N:\; \Lambda^2 T^{1,0}(M)\arrow
T^{0,1}(M)$ the Nijenhuis tensor. We say that
$N$ is {\bf non-degenerate} if $N$ is an isomorphism everywhere.
Then $(M,I)$ is called {\bf an almost complex 6-manifold
with nowhere degenerate Nijenhuis tensor}.

\hfill

\remark
Such manifolds were investigated by R. Bryant.
His results were  presented at a conference
\cite{_Bryant_talk_}, but never published.
The present author  did not attend (unfortunately) 
and was  not aware of his work.

\hfill

The first thing one notices is that the determinant 
$\det N^*$ gives a section 
\[ 
\det N^* \in \Lambda^{3,0}(M)^{\otimes 2}\otimes \Lambda^{3,0}(M)^*.
\] 
Taking 
\begin{equation}\label{_volume_via_N_Equation_} 
   \det N^* \otimes \overline{\det N^*}\in 
   \Lambda^{3,0}(M)\otimes \Lambda^{0,3}(M)=\Lambda^6(M)
\end{equation} 
we obtain a 
nowhere degenerate real volume form $\Vol_I$ on $M$.
This form is called {\bf the canonical volume form associated
with the Nijenhuis tensor}. This gives a functional
$I \stackrel\Psi\arrow \int_M \Vol_I$ on the space of
almost complex structures. One of the purposes of this
paper is to investigate the critical points of the functional
$\Psi$, in the spirit of Hitchin's work 
(\cite{_Hitchin:3-forms_}, \cite{_Hitchin:stable_}).

\subsection{Connections with totally antisymmetric
torsion}
\label{_conne_tota_antisy_Subsection_}

Let $(M, g)$ be a Riemannian manifold, 
$\nabla:\; TM \arrow  TM \otimes \Lambda^1M$ a connection, and
$T \subset \Lambda^2M\otimes TM$ its torsion. Identifying
$TM$ and $\Lambda^1M$ via $g$, we may consider $T$ as an
element in $\Lambda^2M \otimes \Lambda^1 M$, that is, a
3-form on $TM$. If $T$ is totally skew-symmetric as
a  3-form on $TM$, we say that $\nabla$ is {\bf a connection
with totally skew-symmetric} (or {\bf totally antisymmetric})
torsion. If, in addition, $M$ is Hermitian, and
$\nabla$ preserves the Hermitian structure,
we say that $\nabla$ {\bf is a Hermitian 
connection with totally antisymmetric torsion}.

Connections with totally skew-symmetric torsion are
extremely useful in physics and differential geometry.
An important example of such a connection is 
provided by a theorem of Bismut (\cite{_Bismut:connection_}).

\hfill

\theorem
Let $(M,I)$ be a complex manifold, and $g$ a Hermitian
metric. Then $M$ admits a unique connection with totally
skew-symmetric torsion preserving $I$ and $g$.

\endproof

\hfill

Connections with totally skew-symmetric torsion
were studied at great length by Friedrich, Ivanov
and others (see e.g. \cite{_Friedrich_Ivanov:physics_},
\cite{_Friedrich_}, \cite{_Agricola_Friedrich_}).
Bismut's theorem requires the base manifold to be complex.
Motivated by string theory, 
Friedrich and  Ivanov generalized Bismut's theorem
to non-integrable almost complex manifolds
(\cite{_Friedrich_Ivanov:physics_}). For completeness,
we sketch a proof of their theorem below.

\hfill

\theorem\label{_FI_existence_skew_sy_to_Theorem_}
Let $(M, I, \omega)$ be an almost complex Hermitian manifold,
and 
\[
N:\; \Lambda^2 T^{1,0}(M)\arrow T^{0,1}(M).
\]
the Nijenhuis tensor. Consider
the 3-linear form 
$\rho:\; T^{1,0}(M)\times T^{1,0}(M)\times T^{1,0}(M)\arrow \C$,
\begin{equation}\label{_3_for_from_Ni_Equation_}
\rho(x,y,z):= \omega(N(x,y), z)
\end{equation}
Then $M$ admits a
connection $\nabla$ with totally skew-symmetric torsion
preserving $(\omega, I)$ if and only if $\rho$ is skew-symmetric.
Moreover, such a connection is unique.

\hfill

{\bf Sketch of a proof:}
\ref{_FI_existence_skew_sy_to_Theorem_} is proven
essentially
in the same way as one proves Bismut's theorem and existence
and uniqueness of a Levi-Civita connection. Let $(M,I,g)$ 
be a Hermitian manifold, and $\nabla_0$ a Hermitian
connection. Then all Hermitian connections can be
obtained by taking $\nabla(A):= \nabla_0 + A$, where
$A$ is a 1-form with coefficients in the algebra
${\goth u}(TM)$ of all skew-Hermitian endomorpisms.
The torsion $T_A$ of $\nabla(A)$ is written as
\[
T_A= T_0 + \Alt_{12}(A),
\]
where $T_0$ is a torsion of $\nabla_0$, and 
$\Alt_{12}$ denotes the antisymmetrization of
\[ 
   \Lambda^1(M)\otimes {\goth u}(TM)\subset 
   \Lambda^1(M)\otimes\Lambda^1(M)\otimes TM 
\]
over the first two indices. We identify
${\goth u}(TM)$ with $\Lambda^{1,1}(M)$ in
a standard way.
Then \ref{_FI_existence_skew_sy_to_Theorem_} 
can be reinterpreted as a statement about linear-algebraic
properties of the operator
\begin{equation}\label{_Alt_12_on_Herm_conne_Equation_}
\Alt_{12}:\; \Lambda^1(M)\otimes\Lambda^{1,1}(M)\arrow
 \bigg(\Lambda^2(M)\otimes\Lambda^1(M)\bigg)^{(2,1)+(1,2)}.
\end{equation}
(the superscript $(\dots)^{(2,1)+(1,2)}$
means taking ${(2,1)+(1,2)}$-part with respect
to the Hodge decomposition), as follows.

By definition, the Nijenhuis tensor $N$ is a section
of $\Lambda^{2,0}\otimes T^{0,1}$. Identifying
$T^{0,1}$ with $\Lambda^{1,0}$ via $g$, we can
consider $N$ as an element of $\Lambda^{2,0}\otimes
\Lambda^{1,0}$. By Cartan's formula, $N$ is equal
to the $(3,0)$-part of the torsion. Therefore,
existence of a connection with totally skew-symmetric
torsion implies that \eqref{_3_for_from_Ni_Equation_}
is skew-symmetric.

Conversely, assume that \eqref{_3_for_from_Ni_Equation_}
is skew-symmetric. Since
\eqref{_Alt_12_on_Herm_conne_Equation_}
maps $\Lambda^1(M)\otimes\Lambda^{1,1}(M)$ to
$(2,1)\oplus (1,2)$-tensors, the $(3,0)$ and
$(0,3)$-parts of torsion stay skew-symmetric
if we modify the connection by adding
$A\in \Lambda^1\otimes {\goth u}(TM)$. Denote
by $T_1$ the $(2,1)\oplus (1,2)$-part of
the torsion $T_0$.
To prove 
\ref{_FI_existence_skew_sy_to_Theorem_},
we need to find $A \in
\Lambda^1(M)\otimes\Lambda^{1,1}(M)$
such that $T_1-\Alt_{12}(A)$ is totally
skew-symmetric.

The map $\Alt_{12}:\; \Lambda^1(M)\otimes\Lambda^{2}(M)\arrow
 \Lambda^2(M)\otimes\Lambda^1(M)$
is an isomorphism, as a dimension count implies
(this map has no kernel, which is easy to see). 
Therefore, \eqref{_Alt_12_on_Herm_conne_Equation_}
is injective. Using dimension count again,
we find that cokernel of
 \eqref{_Alt_12_on_Herm_conne_Equation_}
projects isomorphically into
\[ \Lambda^{2,1}(M)\oplus\Lambda^{1,2}(M)\subset 
   \Lambda^2(M)\otimes\Lambda^1(M).
\] 
Therefore, for any $T_1$
in $(2,1)\oplus (1,2)$-part of 
$\Lambda^2(M)\otimes\Lambda^1(M)$
there exists $A \in
\Lambda^1(M)\otimes\Lambda^{1,1}(M)$ and
$B\in \Lambda^{2,1}(M)\oplus\Lambda^{1,2}(M)$
such that $T_1=\Alt_{12}(A)+B$.
\endproof

\subsection[Connections with antisymmetric
  torsion on almost complex manifolds]{Connections with antisymmetric
  torsion on almost complex 6-manifolds}

Let $(M,I)$ be an almost complex manifold, $N$ 
its Nijenhuis tensor. To obtain all Hermitian 
connections with totally skew-symmetric 
torsion on $(M,I)$, one needs to find all metrics
$g$ for which the tensor $\omega(N(x,y), z)$ is skew-symmetric.
As \ref{_FI_existence_skew_sy_to_Theorem_} implies, these
metrics are precisely those for which such
a connection exists.

We also prove the following proposition

\hfill

\proposition \label{_exi_metri_Nije_skew_Proposition_}
Let $(M,I)$ be an almost complex 6-manifold with
Nijenhuis tensor which is non-degenerate in a dense
subset of $M$, and $g$ a Hermitian metric
admitting a connection with totally antisymmetric
torsion. Then $g$ is uniquely 
determined by $I$, up to conformal equivalence.
Moreover, the Riemannian metric 
$g$ determines $I$ uniquely, unless
$(M,g)$ is locally isometric to a 6-sphere.

\hfill

{\bf Proof:} This is 
 \ref{_Hermi_tota_antisymme_unique_Proposition_}
and \ref{_NK_determined_by_metric_Proposition_}.
 \endproof

\subsection{Correspondence with the results of R. Bryant}

Since the first version of this paper was written,
the previously unpublished results of R. Bryant appeared 
in a fundamental and important preprint \cite{_Bryant:6dim_}. There is a significant
overlap with our research, although the presentation
and terminology is different. The property 
\eqref{_3_for_from_Ni_Equation_} (which is equivalent
to existence of Hermitian connection with totally 
antisymmetric curvature) is called 
``Nijenhuis tensor of real type'' in \cite{_Bryant:6dim_}.
The main focus of \cite{_Bryant:6dim_}
is the so-called ``quasi-integrable almost complex
manifold'': manifolds with Nijenhuis tensor of real type,
 which is at every point of $M$ either 
non-degenerate (of constant signature) or zero. Examples of 
such structures are found. In particular,
 all twistor spaces of Kaehler surfaces with
sign-definite holomorphic bisectional curvature
are shown to be quasi-integrable.
A variant of \ref{_extre_for_M_I_Theorem_} 
is also proven. It is shown (\cite{_Bryant:6dim_},
Proposition 8) that nearly Kaehler
manifolds are critical points of the functional
$\Vol_I$.


\section{Nearly K\"ahler manifolds: an introduction}
\label{_NK_Intro_Section_}


Nearly K\"ahler manifolds (also known as 
$K$-spaces or almost Tachibana spaces)
were defined and studied by 
Alfred Gray (\cite{_Gray:1965_}, \cite{_Gray:NK_}, 
\cite{_Gray:weak_holo_}, \cite{_Gray:structure_NK_}) in a general
context of intrinsic torsion of $U(n)$-\-struc\-tures
and weak ho\-lo\-no\-mies. An almost complex Hermitian
manifold $(M,I)$ is called {\bf nearly K\"ahler}
if $\nabla_X(I) X=0$, for any vector fields $X$ ($\nabla$ 
denotes the Levi-Civita connection). In other words, the
tensor $\nabla\omega$ must be totally skew-symmetric, for
$\omega$ the Hermitian form on $M$. If
$\nabla_X(\omega)\neq 0$ for any non-zero
vector field $X$, $M$ is called {\bf strictly nearly K\"ahler}.

In this section, we give an overview of known results
and ``folk theorems'' of nearly K\"ahler geometry.
Most of this theory was known (in a different context)
since 1980-ies, when the study of Killing spinors
was initiated (\cite{_Baum_etc_}). 

\subsection{Splitting theorems for nearly K\"ahler manifolds}

As V. F. Kirichenko proved,
nearly K\"ahler manifolds admit a connection with totally
antisymmetric, parallel torsion
(\cite{_Kirichenko:torsion_NK_}). This observation
was used to prove a splitting theorem for 
nearly K\"ahler manifolds: any nearly K\"ahler
manifold is locally a Riemannian
product of a K\"ahler manifold and 
a  strictly nearly K\"ahler one
(\cite{_Gray:structure_NK_}, \cite{_Nagy:NK_}).

A powerful classification theorem for
Riemannian manifolds admitting an orthogonal 
connection with irreducible connection and parallel
torsion was obtained by R. Cleyton and A. Swann 
in \cite{_Cleyton_Swann_}.
Cleyton and Swann proved that any such 
manifold is either locally homogeneous, has vanishing torsion,
or has weak holonomy $G_2$ (in dimension 7) or
$SU(3)$ (in dimension 6).

Using Kirichenko theorem, this result can be used to obtain
a classification of nearly K\"ahler manifolds.
P.-A. Nagy (\cite{_Nagy:splittine_}) has shown that 
that any strictly nearly K\"ahler manifold
is locally a product of locally homogeneous manifolds,
strictly nearly K\"ahler 6-manifolds, and twistor
spaces of quaternionic K\"ahler manifolds of
positive Ricci curvature, equipped with
the Eells-Salamon metric.

These days the term ``nearly K\"ahler'' usually denotes
strictly nearly K\"ahler 6-manifolds. In sequel we shall follow
this usage, often omitting ``strictly'' and ``6-dimensional''.

In dimension 6, a manifold is (strictly) nearly 
K\"ahler if and only if it admits a Killing spinor
(\cite{_Grunewald_}). Therefore, such a manifold is
Einstein, with positive Einstein constant. 

As one can easily show (see \ref{_NK_from_forms_Theorem_}), 
strictly nearly K\"ahler 6-manifolds can be defined as
6-manifolds with structure group $SU(3)$ and fundamental
forms $\omega\in \Lambda^{1,1}_\R(M)$, 
$\Omega\in \Lambda^{3,0}(M)$, satisfying 
$d \omega=3\lambda\Re\Omega$,
$d \Im\Omega = -2 \lambda\omega^2$.
 An excellent
introduction to nearly K\"ahler geometry is
found in \cite{_Moroianu_Nagy_Semmelmann:NK_}.

The most puzzling aspect of nearly K\"ahler geometry
is a complete lack of non-homogeneous 
examples. With the exception of 4 homogeneous 
cases described below (Subsection
\ref{_exa_NK_intro_Subsection_}), no other compact examples of
strictly nearly K\"ahler 6-manifolds are known to exist.

\subsection{Nearly K\"ahler manifolds 
in $G_2$-geometry and physics}
\label{_G_2_conical_Subsection_}

Nearly K\"ahler manifolds have many uses in geometry
and physics. Along with Calabi-Yau manifolds, 
nearly K\"ahler manifolds appear as  target spaces
for supersymmetric sigma-models, solving equations
of type II string theory. These manifolds are the only
6-manifolds admitting a Killing spinor. This implies
that a Riemannian cone $C(M)$  of a nearly K\"ahler manifold
has a parallel spinor. 

Let $(M, g)$ be a Riemannian manifold. Recall
that the {\bf Riemannian cone} of $(M, g)$ is a product
$M\times \R^{> 0}$, with a metric 
$gt^2 \oplus \lambda \cdot dt^2$,
where $t$ is a unit parameter on $\R^{>0}$,
and $\lambda$ a constant.
It is well known that $M$ admits a real Killing spinor
if and only if $C(M)$ admits a parallel spinor
(for appropriate choice of $\lambda$).
Then, $C(M)$ has restricted holonomy, for any nearly K\"ahler
6-manifold. It is easy to check that in fact $C(M)$
has holonomy $G_2$. This explains a tremendous importance that 
nearly K\"ahler manifolds play in $G_2$-geometry.

We give a brief introduction of $G_2$-geometry, following
\cite{_Hitchin:3-forms_} and \cite{_Joyce_Book_}.
Let $V^7$ be a 7-dimensional real vector space.
The group $GL(7, \R)$ acts on $\Lambda^3(V^7)$
with two open orbits. For $\nu$ in one of these orbits,
its stabilizer $St(\nu)\subset GL(7, \R)$ is 14-dimensional,
as a dimension count insures. It is easy to check that
$St(\nu)$ is a real form of a Lie group $G_2$. 
For one of these orbits, $St(\nu)$ is a compact
form of  $G_2$, for another one it is non-compact.
A 3-form $\nu\in \Lambda^3(V^7)$ is called 
{\bf stable} if its stabilizer is 
a compact form of  $G_2$. 

A 7-manifold $X$ equipped with a 3-form $\rho$
is called {\bf a $G_2$-manifold} if $\rho$ is stable
everywhere in $X$. In this case, the structure
group of $X$ is reduced to $G_2$. Also, $X$ is
equipped with a natural Riemannian structure:
\begin{equation}\label{_metric_G_2_Equation_}
x, y \arrow 
   \int_X (\rho\cntrct x) \wedge (\rho\cntrct x)\wedge
   \rho \ \  (x, y \in TM).
\end{equation}
A $G_2$-manifold is called {\bf parallel} if
$\nabla\rho=0$, where $\nabla$ is the Levi-Civita
connection associated with this Riemannian structure.

Isolated singularities of $G_2$-manifolds are of paramount
importance in physics (\cite{_Acharya_Gukov_}, \cite{_Atiyah_Witten_}).
A simplest example of an isolated singular point is a conical singularity. 

A metric space $X$ with marked points 
$x_1, ... x_n$ is called {\bf a space with isolated
  singularities}, if $X\backslash \{x_1, ... x_n\}$
is a Riemannian manifold. Consider a space $(X,x)$
with a single singular point. The singularity
$x\in X$ is called {\bf conical} if $X$
is equipped with a flow acting on $X$ by 
homotheties and 
contracting $X$ to $x$. In this case,
$X\backslash x$ is isomorphic to a Riemannian
cone of a Riemannian manifold $M$.

It is easy to check that the cone $C(M)$ of
a nearly K\"ahler manifold is equipped with a parallel
$G_2$-structure, and, conversely, every conical
singularity of a parallel $G_2$-manifold
is obtained as $C(M)$, for some  nearly K\"ahler manifold
$M$ (\cite{_Hitchin:stable_}, \cite{_Ivanov_Parton_Piccini_}). 
For completeness' sake, we give a sketch of a 
proof of this result in \ref{_G_2_on_cone_Proposition_}.

The idea of this correspondence is quite clear.
Let $X=C(M)$ be a parallel $G_2$-manifold,
and $\omega_C$ its 3-form. Unless $X$ is flat,
we may assume that
$X$ has holonomy which is equal to $G_2$ and not
its proper subgroup. Indeed, if holonomy of $X$ is less
than $G_2$, by Berger's classification of irreducible
holonomies $X$ is represented (as a Riemannian manifold)
as a product of manifolds of smalled dimension.
However, the singular point of the metric 
completion $\bar X$ is isolated,
and this precludes such a decomposition, unless
$\bar X$ is smooth. In the latter case,
$X$ is flat.

Since holonomy of $X$ is (strictly) $G_2$,
the 3-form can be reconstructed from the
Riemannian structure uniquely. After rescaling,
we may assume that the Riemannian structure
structure on $X=C(M)$ is homogeneous of weight 2,
with respect to the action of $\R^{>0}$ on $C(M)$.
Then $\omega_C$ is homogeneous of weight 3.
Homogeneous $G_2$-structures on $C(M)$ correspond
naturally to $SU(3)$-structures on $M$. 
We write $\omega_C$ as $t^2 \pi^*\omega \wedge dt + t^3 \pi^*\rho$,
where $\rho$, $\omega$ are forms on $M$, and 
$\pi:\; C(M) \arrow M$ is the standard projection.
From a local coordinate expression of a $G_2$-form,
we find that $\omega$ is a Hermitian form
corresponding to an almost complex structure $I$,
and $\rho=\Re\Omega$ for a nowhere degenerate
$(3,0)$-form $\Omega$ on $(M,I)$. 

The converse is proven by the same computation: given 
an $SU(3)$-manifold $(M,I, \omega, \Omega)$,
we write a 3-form 
\begin{equation}\label{_omega_C_on_cone_Equation_} 
  \omega_C:= t^2 \pi^*\omega \wedge dt + t^3 \pi^*\rho,
\end{equation}
on $C(M)$, and show that it is a $G_2$-structure,
using a coordinate expression for a $G_2$-form.

As Fernandez and Gray proved in 
\cite{_Fernandez_Gray:G_2_}, a $G_2$-manifold
$(X, \omega_C)$ is parallel if and only if $\omega_C$
is harmonic. For the form
\eqref{_omega_C_on_cone_Equation_},
$d\omega_C=0$ is translated into
\begin{equation}\label{_d_omega_NK_Equation_}
  d\omega=3\rho.
\end{equation} 

Since $*\rho = I\rho$ and 
$*\omega = \omega^2$, the condition $d^*\omega_C=0$
becomes $d I \rho = -2 \omega^2$.
After an appropriate rescaling,
we find that this is precisely the condition
defining the nearly K\"ahler structure 
(\ref{_NK_from_forms_Theorem_}).
Therefore,
$C(M)$ is a $G_2$-manifold if and only if $M$
is nearly K\"ahler (see Subsection 
\ref{_G_2_on_cones_Subsection_} for a 
more detailed argument).

The correspondence between conical singularities 
of $G_2$-manifolds and nearly K\"ahler geometry
can be used further to study the locally conformally
parallel $G_2$-manifolds (see also \cite{_Ivanov_Parton_Piccini_}).
{\bf Locally conformally parallel $G_2$-manifold} is a
7-manifold $M$ with a covering $\tilde M$ equipped
with a parallel $G_2$-structure, with the deck
transform acting on $\tilde M$ by homotheties.
Since homotheties preserve the Levi-Civita connection
$\tilde \nabla$ on $\tilde M$, $\tilde \nabla$
descends to a torsion-free connection on $M$,
which is no longer orthogonal, but preserves the
conformal class of a metric. Such a connection
is called {\bf a Weyl connection}, and a conformal
manifold of dimension $>2$ equipped with a torsion-free connection
preserving the conformal class {\bf a Weyl manifold}.
The Weyl manifolds are a subject of much study in
conformal geometry (see e.g. \cite{_Dragomir_Ornea_}
and the reference therein).

The key theorem of Weyl geometry is proven by
P. Gauduchon (\cite{_Gauduchon_1984_}). He has shown 
that any compact Weyl manifold is equipped with a privileged
metric in its conformal class. This metric
(called  {\bf a Gauduchon metric} now) is defined as
follows.

Let $(M, [g], \nabla)$ be a compact Weyl manifold, where
$[g]$ is a conformal class, and $g\in [g]$ any metric
within this conformal class. Since $\nabla[g]=0$,
we have $\nabla(g)= g \otimes \theta$, where $\theta$
is a 1-form, called {\bf a Lee form}. A metric
$g$ is called {\bf Gauduchon} if $\theta$ satisfies
$d^*\theta=0$. A Gauduchon metric is unique (up to a
complex multiplier). 

Let now $(M, \nabla, [g])$ 
be a Weyl manifold with a Ricci-flat connection $\nabla$.
In \cite{_Gauduchon:Weyl_Eisntein_}, 
Gauduchon has shown that the Lee form
$\theta$ of the Gauduchon metric on $M$
is parallel with respect to the Levi-Civita connection
associated with this metric. 

Applying this argument to a compact
locally conformally parallel $G_2$-manifold $M$, we
obtain that the Lee form is parallel.
From this one infers that 
the parallel $G_2$-covering $\tilde M$ 
of $M$ is a cone over some Riemannian manifold $S$
(see e.g. \cite{_Verbitsky:LCHK_}, Proposition 11.1,
also see \cite{_Kamishima_Ornea_} and 
\cite{_Gini_Ornea_Parton_}). 
Using the argument stated above, we find that 
this manifold is in fact nearly K\"ahler.
Therefore, $S$ is Einstein, with positive
Ricci curvature. Since $M$ is compact,
$S$ is complete, and by Myers theorem,
$S$ is actually compact (see \cite{_Verbitsky:LCHK_},
Remark 10.7). Now, the argument which proves
Theorem 12.1 of \cite{_Verbitsky:LCHK_} can be used
to show that $\dim H^1(M, \Q)=1$, and 
$M = C(S)/\Z$. This gives the following
structure theorem, which is proven independently in
\cite{_Ivanov_Parton_Piccini_}.

\hfill

\theorem
Let $M$ be a compact locally conformally parallel
$G_2$-\-ma\-ni\-fold. Then $M= C(S)/\Z$, where
$S$ is a nearly K\"ahler manifold, and the
$\Z$-action on $C(S)\cong S \times \R^{>0}$ is generated by
a map $(x, t) \arrow (\phi(x), qt)$, where
$|q|>1$ is a real number, and $\phi:\; S \arrow S$
an automorphism of nearly K\"ahler structure.

\endproof

\subsection{Examples of nearly K\"ahler manifolds}
\label{_exa_NK_intro_Subsection_}

Just as the conical singularities of parallel 
$G_2$-manifolds correspond to nearly K\"ahler manifolds,
the conical singularities of $Spin(7)$-\-ma\-ni\-folds
correspond to the so-called ``nearly parallel'' 
$G_2$-manifolds (see \cite{_Ivanov:torsion_G_2_}).
A $G_2$-manifold $(M, \omega)$ is called {\bf nearly parallel}
if $d\omega = c *\omega$, where $c$ is some constant.
The analogy between nearly K\"ahler 6-manifolds and nearly
parallel $G_2$-manifolds is almost perfect.
These manifolds admit a connection with totally
antisymmetric torsion and have weak holonomy
$SU(3)$ and $G_2$ respectively.
N. Hitchin realized nearly K\"ahler 6-manifolds and nearly
parallel $G_2$-manifolds as extrema of a certain 
functional, called {\bf Hitchin functional} by physicists
(see \cite{_Hitchin:stable_}). 

However, examples of nearly parallel $G_2$-manifolds
are found in profusion (every 3-Sasakian manifold
is nearly parallel $G_2$), and compact nearly K\"ahler
manifolds are rare.

Only 4 compact examples are known
(see the list below); 
all of them homogeneous.
In \cite{_Butruille_} it was shown that any 
homogeneous nearly K\"ahler 6-manifold 
belongs to this list. 

\begin{enumerate}

\item The 6-dimensional sphere $S^6$. 
Since the cone  $C(S^6)$ is flat, $S^6$ is a nearly 
K\"ahler manifold, as shown in Subsection
\ref{_G_2_conical_Subsection_}. The 
almost complex structure on $S^6$ is reconstructed
from the octonion action, and the metric is standard.

\item $S^3\times S^3$, with the complex structure
mapping $\xi_i$ to $\xi'_i$,  $\xi_i'$ to $-\xi_i$,
where $\xi_i$, $\xi'_i$, $i=1,2,3$ is a basis
of left invariant 1-forms on the first and the second 
component. 

\item Given a self-dual Einstein Riemannian 4-manifold $M$
with positive Einstein constant, one defines
its {\bf twistor space} $\Tw(M)$ as a total space of
a bundle of unit spheres in $\Lambda^2_-(M)$
of anti-self-dual 2-forms. Then $\Tw(M)$ 
has a natural K\"ahler-Einstein structure $(I_+, g)$,
obtained by interpreting unit vectors
in  $\Lambda^2_-(M)$ as  complex structure
operators on $TM$. Changing the sign of 
$I_+$ on $TM$, we obtain an almost complex structure $I_-$
which is also compatible with the metric $g$
(\cite{_Eells_Salamon_}). 
A straightforward computation insures that
$(\Tw(M), I_-, g)$ is nearly K\"ahler
(\cite{_Muskarov_}).

As N. Hitchin proved (\cite{_Hitchin:kah_twistors_}),
there are only two compact self-dual
Einstein 4-manifolds: $S^4$ and $\C P^2$.
The corresponding twistor spaces are
$\C P^3$ and the flag space $F(1,2)$.
The almost complex structure operator $I_-$ 
induces a nearly K\"ahler structure on these
two symmetric spaces. 

\end{enumerate}

\subsection[Nearly K\"ahler manifolds as extrema of
intrinsic  volume] {Nearly K\"ahler manifolds are extrema of volume on
almost complex manifolds with nowhere degenerate Nijenhuis
tensor}

Let $(M, I, \omega)$ be a nearly K\"ahler manifold, and
$N^*:\; \Lambda^{0,1}(M) \arrow \Lambda^{2,0}(M)$
the Nijenhuis tensor. By Cartan's formula, $N^*$ is the
$(2,-1)$-part of the de Rham differential (with 
respect to the Hodge decomposition). 
In \ref{_NK_from_forms_Theorem_} it is shown that 
$d \omega$ is a real part of a nowhere degenerate
$(3,0)$-form $\Omega$. Therefore, the 3-form
\[
\omega(N(x,y), z) = d\omega(x,y,z) =  \Re\Omega(x,y,z)
\]
is nowhere degenerate on $T^{1,0}(M)$.
We obtain that the Nijenhuis tensor 
$N$ is nowhere degenerate.

The main result of this paper is the following theorem,
which is analogous to \cite{_Hitchin:stable_}.

\hfill

\theorem \label{_extre_for_M_I_Theorem_}
Let $(M,I)$ be a compact almost complex 
6-manifold with nowhere degenerate Nijenhuis tensor
admitting a Hermitian connection with totally
antisymmetric torsion.
Consider the functional 
\begin{equation}\label{_volume_func_in_intro_theorem_Equation_}
I\arrow \int_M\Vol_I
\end{equation}
on the space of such manifolds constructed in Subsection
\ref{_Nije_non-dege_Subsection_}.
Then \eqref{_volume_func_in_intro_theorem_Equation_}
has a critical point at $I$ if and only if $(M,I)$
admits a nearly K\"ahler metric.

\hfill

{\bf Proof:} Follows from 
\ref{_extre_of_Psi_via_d_omega_Proposition_}
and \ref{_NK_from_forms_Theorem_}. \endproof

\hfill

\remark
As follows from \ref{_NK_unique_from_I_Corollary_},
the nearly K\"ahler metric on $(M,I)$ is uniquely determined by
the almost complex structure.


\section[Almost complex structures on 6-manifolds 
and connections with torsion]{Almost complex 
structures and connections \\
with totally antisymmetric torsion}
\label{_conne_antisymme_Nije_Section_}


Let $(M,I)$ be a 6-dimensional almost complex manifold,
and \[ N^*:\; \Lambda^{0,1}(M) \arrow \Lambda^{2,0}(M)\]
its Nijenhuis tensor. Given a point $x\in M$, the operator
$N^*\restrict{\Lambda^{0,1}_x(M)}$ can a priori take any
value within $\Hom(\Lambda^{0,1}(M), \Lambda^{2,0}(M))$.
For $N^*\restrict{\Lambda^{0,1}_x(M)}$
generic, the stabilizer $St(N^*_x)$ of $N^*_x$
within $GL(T_xM)$ is 2-dimensional. If we fix a 
complex parameter, the eigenspaces of
$N^*_x$ (taken in apporpriate sense)
define a frame in $TM$. Thus, a geometry
of a ``very generic'' 6-dimensional almost complex manifold
is rather trivial.

However, for a $N^*_x$ inside a 10-dimensional subspace
\[ W_0\subset \Hom(\Lambda^{0,1}(M), \Lambda^{2,0}(M)),\]
(\ref{_W_0_Remark_}), the stabilizer
$St(N^*_x)$ contains $SU(3)$, and the geometry of $(M,I)$ becomes
more interesting. 

\hfill

\proposition\label{_Hermi_tota_antisymme_unique_Proposition_}
Let $(M,I)$ be an almost complex 6-manifold with 
Nijenhuis tensor which is non-degenerate 
in a dense set. Assume that $(M,I)$ admits a Hermitian
structure $\omega$ and a Hermitian connection
with totally antisymmetric torsion. Then
$\omega$ is uniquely determined by $I$, up
to conformal equivalence.

\hfill

{\bf Proof:}
Consider the map 
\begin{equation}\label{_C_defi_Equation_}
   C:= \Id\otimes N^*:\; \Lambda^{1,1}(M)\arrow  
   \Lambda^{1,0}(M)\otimes \Lambda^{2,0}(M)
\end{equation}
obtained by acting with the Nijenhuis tensor
$N^*:\; \Lambda^{0,1}(M) \arrow \Lambda^{2,0}(M)$
on the second tensor multiplier of
$\Lambda^{1,1}(M)\cong \Lambda^{1,0}(M)\otimes \Lambda^{0,1}(M)$.
Then $C$ maps $\omega$ to a 3-form 
\[ x, y, z \arrow \omega(N(x,y), z).\]
As \ref{_FI_existence_skew_sy_to_Theorem_} implies,
$(M, I, \omega)$ admits a Hermitian connection
with totally antisymmetric torsion if and only if
$C(\omega)$ lies inside a 1-dimensional space
\[ \Lambda^{3,0}(M)\subset \Lambda^{1,0}(M)\otimes \Lambda^{2,0}(M).
\]
However, $C$ is an isomorphism in a dense subset of $M$,
hence, all $\omega$ which satisfy the conditions
of \ref{_FI_existence_skew_sy_to_Theorem_} are proportional.
\endproof

\hfill

\remark\label{_W_0_Remark_}
The same argument proves that
an almost complex manifold admits
a Hermitian connection
with totally antisymmetric torsion
if and only if $C^{-1}(\Lambda^{3,0}(M))$
contains a Hermitian form. This is
the space $W_0$ alluded to in the beginning
of this section.

\hfill

\ref{_Hermi_tota_antisymme_unique_Proposition_}
leads to the following corollary.

\hfill

\corollary\label{_NK_unique_from_I_Corollary_} 
Let $(M,I)$ be an almost complex 6-manifold.
Then $(M,I)$ admits at most one strictly nearly
K\"ahler metric, up to a constant multiplier.

\hfill

{\bf Proof:}
Let $\omega_1$ and $ \omega_2$
be nearly K\"ahler metrics on $(M,I)$.
Since $(M,I,\omega_i)$ is strictly nearly K\"ahler,
the 3-form $C(\omega_i)\in \Lambda^{3,0}(M)$
is nowhere degenerate (see \eqref{_C_defi_Equation_}).
Therefore, $(M,I)$ has nowhere degenerate Nijenhuis
tensor. Then, by \ref{_Hermi_tota_antisymme_unique_Proposition_},
$\omega_i$ are proportional: $\omega_1 = f\omega_2$.
However, $d\omega_i^2=0$ on any nearly K\"ahler
3-manifold  (see e.g. \ref{_NK_from_forms_Theorem_} (ii)).
Then $2fdf \wedge \omega_2^2=0$. This implies
$df=0$, because the map $\eta \arrow \eta \wedge \omega_2^2$
is an isomorphism on $\Lambda^1(M)$.
\endproof

\hfill

\remark
The converse is also true:
unless $(M,g)$ is locally isometric to a 6-sphere,
the Riemannian metric $g$ determines the nearly K\"ahler
almost complex structure $I$ uniquely
(\ref{_NK_determined_by_metric_Proposition_}).


\section[Nearly K\"ahler geometry and connections
with torsion]{Nearly K\"ahler geometry and Hermitian connections with 
totally antisymmetric torsion}
\label{_NK_and_G2_antisy_Section_}


\subsection{Hermitian structure on $\Lambda^{3,0}(M)$
and nearly K\"ahler manifolds}

Let $(M,I)$ be an almost complex 6-manifold,
and $\Omega \in \Lambda^{3,0}(M)$ a non-degenerate $(3,0)$-form.
Then $\Omega\wedge \bar\Omega$ is a positive volume
form on $M$. This gives a $\Vol(M)$-valued Hermitian
structure on $\Lambda^{3,0}(M)$. 
If $M$ is in addition Hermitian,
then $M$ is equipped with a natural volume form $\Vol_h$
associated with the metric, and the map
\[ 
\Omega \arrow \frac{\Omega\wedge \bar\Omega}{\Vol_h}
\]
can be considered as a Hermitian metric on
$\Lambda^{3,0}(M)$.
This metric agrees with the usual Riemann-Hodge pairing
known from algebraic geometry, when $I$ is integrable. 
The following definition is a restatement of the classical 
one, see Subsection \ref{_NK_Intro_Section_}.

\hfill

\definition\label{_NK_tens_Definition_}
Let $(M, I, \omega)$ be an almost complex 
Hermitian manifold, and $\nabla$ the Levi-Civita 
connection. Then $(M, I, \omega)$
is called {\bf nearly K\"ahler} if the tensor
$\nabla\omega$ is totally antisymmetric:
\[ \nabla\omega\subset \Lambda^3(M).
\]

\hfill

The following theorem is a main result of this section.

\hfill

\theorem\label{_NK_from_forms_Theorem_}
Let $(M,I, \omega)$ be an almost complex  Hermitian
6-manifold equipped with a $(3,0)$-form $\Omega$. Assume
that $\Omega$ satisfies
$3\lambda\Re\Omega = d\omega$, and $|\Omega|_\omega=1$,
where $\lambda$ is a constant, and $|\cdot|_\omega$
is the Hermitian metric on $\Lambda^{3,0}(M)$
constructed above. Then the following conditions
are equivalent. 
\begin{description}

\item[(i)] $M$ admits a Hermitian 
connection with totally antisymmetric torsion.

\item[(ii)] $d \Omega= -2 \1\lambda\omega^2$

\item[(iii)] $(M,I, \omega)$
is nearly  K\"ahler, and 
$d\omega=\nabla\omega$.
\end{description}

\hfill

The equivalence of (ii) and (iii) is known (see e.g.
\cite{_Hitchin:stable_}, the second part of the proof
of Theorem 6). 

The existence of Hermitian 
connections with totally antisymmetric torsion
on nearly K\"ahler manifolds is also well known
(see Section \ref{_NK_Intro_Section_}). This connection
is written as $\nabla_{NK}=\nabla+T$, where
$\nabla$ is the Levi-Civita connection 
on $M$, and $T$ the operator obtained from
the 3-form $3\lambda\Im\Omega$ by raising one
of the indices. The torsion of  $\nabla_{NK}$ is totally antisymmetric
by construction (it is equal $T$). Also by construction,
we find that $T(\omega)= -3\lambda\Re\Omega$,
hence $\nabla_{NK}(\omega)=0$. Therefore,
$\nabla_{NK}$ is a Hermitian connection 
with totally antisymmetric torsion.
This takes care of the implication
(iii) $\Rightarrow$ (i). 

To prove
\ref{_NK_from_forms_Theorem_}, 
it remains to prove that
(i) implies (ii);
we do that in Subsection 
\ref{_conne_Nije_d_Omega_Subsection_}.
For completeness' sake, we sketch 
the proof of the implication
(ii) $\Rightarrow$ (iii)  in
Subsection \ref{_Nearly_Kah_G_2_Subsection_}.

\hfill

\remark
As \ref{_NK_unique_from_I_Corollary_} 
 shows, a non-K\"ahler
nearly K\"ahler metric on $M$ is uniquely
determined by the almost complex structure $I$.

\subsection{Connections with totally antisymmetric torsion
and Nijenhuis tensor}
\label{_conne_Nije_d_Omega_Subsection_}

\ \ \lemma\label{_d_Omega_omega^2_Lemma_}
In assumptions of \ref{_NK_from_forms_Theorem_},
(i) implies (ii).

\hfill

{\bf Proof.}

\hfill

{\bf Step 1:} {\bf We show that $d\Omega \in \Lambda^{2,2}(M)$.}

Were $(M,I)$ integrable, the differential $d$ would have
only (0,1)- and (1,0)-part with respect to the Hodge
decomposition: $d=d^{1,0}+d^{0,1}$. For a general almost complex manifold,
$d$ splits onto 4 parts:
\[
d = d^{2,-1}+d^{1,0}+d^{0,1}+d^{-1,2}.
\]
This follows immediately from the Leibniz rule.
However,
\begin{equation}\label{_d^2_Omega_vanishes_Equation_}
0 = d^2 \omega = d(\Omega+ \bar \Omega) = d \Omega + d \bar\Omega.
\end{equation}
Since $\Lambda^{p,q}(M)$ vanishes for $p$ or $q>3$,
we also have
\begin{equation}\label{_4_terms_d_Omega_Equation_}
d \Omega + d \bar\Omega= d^{0,1}\Omega + d^{-1,2}\Omega + 
  d^{2,-1}\bar\Omega + d^{1,0}\bar\Omega
\end{equation}
The four terms on the right hand side of 
\eqref{_4_terms_d_Omega_Equation_} have Hodge types
$(3,1)$, $(2,2)$, $(2,2)$ and $(1,3)$. Since their
sum vanishes by \eqref{_d^2_Omega_vanishes_Equation_},
we obtain
\[
d^{0,1}\Omega =0, \ \ d^{1,0}\bar\Omega=0, \ \ 
d^{2,-1}\bar\Omega=-d^{-1,2}\Omega.
\]
Then \eqref{_4_terms_d_Omega_Equation_}
gives
\begin{equation}\label{_d_Omega_via_d^-1,2_Equation_}
d\Omega =-d^{2,-1}\bar\Omega=d^{-1,2}\Omega.
\end{equation}

{\bf Step 2:} 
\begin{equation}\label{_d^2,-1_via_Nije_Equation_}
d^{2,-1}\restrict{\Lambda^{1,1}(M)}=\bigwedge\circ\Id\otimes N^*,
\end{equation}
where $N^*:\; \Lambda^{0,1}(M)\arrow \Lambda^{2,0}(M)$
is the Nijenhuis tensor, 
\[ 
  \Id\otimes N^*:\; \Lambda^{1,1}(M)\arrow 
  \Lambda^{2,0}(M)\otimes \Lambda^{1,0}(M)
\] 
acts as $N^*$ on the second multiplier of 
$\Lambda^{1,1}(M)\cong \Lambda^{1,0}(M)\otimes \Lambda^{0,1}(M)$,
and $\bigwedge$ denotes the exterior product.
\eqref{_d^2,-1_via_Nije_Equation_} is immediately
implied by the Cartan's formula for the de Rham differential.

\hfill

{\bf Step 3:} From the existence of Hermitian connection
with totally antisymmetric torsion we obtain that the form
\[ \omega(N(x,y), z):\; T^{1,0}M\times T^{1,0}M\times T^{1,0}M\arrow \C\]
is totally antisymmetric (see \ref{_FI_existence_skew_sy_to_Theorem_}).
From  \eqref{_d^2,-1_via_Nije_Equation_} it follows that
\begin{equation}\label{_N_via_omega_Omega_Equation_}
\omega(N(x,y), z)= d\omega=3\lambda \Re\Omega.
\end{equation}
Consider an orthonormal frame $d z_1, d z_2, d z_3$
in $\Lambda^{1,0}(M)$, satisfying $\Omega= dz_1\wedge dz_2 \wedge dz_3$
(such a frame exists because $|\Omega|_\omega=1$).
Then \eqref{_N_via_omega_Omega_Equation_} gives
\begin{equation}\label{_N_basis_Equation_}
N^*(d\bar z_i)= \lambda d \check z_i,
\end{equation}
where
$d\check z_1= dz_2 \wedge dz_3$, $d\check z_2= -dz_1 \wedge dz_3$, 
$d\check z_3= dz_1 \wedge dz_2$. 

\hfill

{\bf Step 4:} 
Using Cartan's formula as in Step 2, we express
$d^{-1,2}\Omega$ through the Nijenhuis tensor.
Then \eqref{_d_Omega_via_d^-1,2_Equation_}
can be used to write $d\Omega =d^{-1,2}\Omega$
in terms of $N^*$. Finally, \eqref{_N_basis_Equation_} 
allows to write $d^{-1,2}\Omega$ in coordinates,
obtaining $d\Omega = -2\1 \lambda \omega^2$.
\endproof

\subsection{$G_2$-structures on cones of 
Hermitian 6-manifolds}
\label{_G_2_on_cones_Subsection_}

\proposition\label{_G_2_on_cone_Proposition_}
Let $(M, I, \omega)$ be an almost complex
Hermitian manifold, $\Omega\in \Lambda^{3,0}(M)$
a $(3,0)$-form which satisfies $d\omega =
3\lambda\Re\Omega$, for some real 
constant, and $|\Omega|_\omega=1$.
Assume, in addition, that $d\Omega =-2\1\lambda\omega^2$.
Consider the cone $C(M)= M\times \R^{>0}$, equipped
with a 3-form $\rho= 3 t^2 \omega \wedge dt +  t^3 d\omega$,
where $t$ is the unit parameter on the $\R^{>0}$-component.
Then $(C(M), \rho)$ is a parallel, $G_2$-manifold
(see Subsection \ref{_G_2_conical_Subsection_}).
Moreover, any parallel $G_2$-structure $\rho'$ on
$C(M)$ is obtained this way, assuming that $\rho'$
is homogeneous of weight 3 with respect to the
the natural action of $\R^{>0}$ on $C(M)$.

\hfill

{\bf Proof:} As Fernandez and Gray has shown
(\cite{_Fernandez_Gray:G_2_}), to show that
a $G_2$-structure $\rho$ is parallel it suffices
to prove that $d\rho = d^*\rho =0$. Clearly,
$d\rho=0$, because
\[ d \rho = 3  t^2 d \omega \wedge dt + 3  t^2 dt \wedge d\omega=0.
\]
On the other hand, $*(\omega\wedge dt) = \frac 1 2 t^2 \omega^2$, 
and $* d\omega = -3dt \wedge I(d\omega)$, 
where $*$ is taken with respect to the cone metric on $C(M)$.
This is clear, because
$(\omega, \Omega)$ defines an $SU(3)$-structure on $M$,
and $d\omega = 3\lambda \Re\Omega$. Then
\begin{equation}\label{_*_rho_expli_Equation_}
*\rho = \frac 3 2 t^4 \omega^2 -3 t^3 dt \wedge I(d\omega).
\end{equation}
Since $d\Omega =-2\1\lambda\omega^2$ and
$3\lambda d \Re\Omega = d^2 \omega=0$,
we obtain $d \Im\Omega = -2 \lambda\omega^2$.
This gives $d I(d\omega)=-2 \omega^2$,
because $\lambda I(d\omega)=\Im\Omega$.
Then \eqref{_*_rho_expli_Equation_}
implies
\[
d(*\rho)= 6t^3 dt \wedge \omega^2 +3 t^3 dt \wedge d I(d\omega)
 = 6t^3 dt \wedge \omega^2-6 t^3 dt \wedge \omega^2 =0.
\]
We proved that $C(M)$ is a parallel $G_2$-manifold.
The converse statement is straightforward.
\endproof

\hfill

In Subsection \ref{_G_2_conical_Subsection_}
it is shown that the holonomy of $C(M)$ is 
strictly $G_2$, unless it is flat (in the latter case,
$M$ is locally isometric to a sphere). Therefore, 
\ref{_G_2_on_cone_Proposition_} implies the following
corollary.

\hfill

\corollary\label{_G_2_unique_Corollary_}
In assumptions of \ref{_G_2_on_cone_Proposition_},
the almost complex structure is uniquely determined by
the metric, unless $M$ is locally isometric to a 6-sphere.

\endproof

\subsection{Near K\"ahlerness obtained from $G_2$-geometry}
\label{_Nearly_Kah_G_2_Subsection_}

Now we can conclude the proof of \ref{_NK_from_forms_Theorem_},
implying \ref{_NK_from_forms_Theorem_} (iii)
from \ref{_NK_from_forms_Theorem_} (ii).
Let $M$ be a 6-manifold satisfying assumptions of 
\ref{_NK_from_forms_Theorem_} (ii). 
Consider the cone $C(M)$ equipped with a parallel
$G_2$-structure $\rho$ as in \ref{_G_2_on_cone_Proposition_}.
Let $g_0$ be a cone metric on $C(M)$. From the argument
used to prove \ref{_G_2_on_cone_Proposition_}, it is 
clear that $g_0$ is a metric induced by the 3-form $\rho$
as in \eqref{_metric_G_2_Equation_}. 

Consider the map
$C(M)\stackrel \tau \arrow M\times \R$ induced by $(m, t) \arrow (m, \log t)$,
and let $g_1=\tau^* g_\pi$ be induced by the product 
metric $g_\pi$ on $M\times \R$. Denote by $\nabla_0$,
$\nabla_1$ the corresponding Levi-Civita connections.
We know that $\nabla_0(\rho)=0$, and we need
to show that 
\begin{equation}\label{_product_conne_to_omega_Equation_}
\nabla_1(\omega)=d\omega.
\end{equation}

The metrics $g_0$, $g_1$ are proportional:
$g_1 = g_0 e^{-t}$. This allows one to relate the
Levi-Civita connections $\nabla_1$ and $\nabla_0$
(see e.g. \cite{_Ornea:LCHK_}):
\[
\nabla_1 = \nabla_0 + \frac1 2 A,
\]
where $A:\; TM \arrow \End(\Lambda^1(M))$
is an $\End(\Lambda^1(M))$-valued 1-form mapping
$X\in TM$ to
\begin{equation}\label{_conne_confo_Equation_}
(\theta, X)\Id - X\otimes \theta + X^\sharp\otimes \theta^\sharp
\end{equation}
and $\theta$ the 1-form defined by $\nabla_0(g_1)=g_1 \otimes \theta$,
 $X\otimes \theta$ the tensor product of $X$ and $\theta$
considered as an endomorphism of $\Lambda^1(M)$, and
$X^\sharp\otimes \theta^\sharp$ the dual endomorphism.

From \eqref{_conne_confo_Equation_} and $\nabla_0(\rho)=0$
we obtain
\begin{equation}\label{_nabla_1_of_rho_Equation_}
(\nabla_1)_X(\rho) = (X, \theta) \rho - (\rho\cntrct X) \wedge \theta
+(\rho \cntrct {\theta^\sharp})\wedge X^\sharp.
\end{equation}
Since $\theta= \frac{dt}{t}$, we have 
$\nabla_1(\theta)=0$, and 
$\nabla_1$ preserves the decomposition
$\Lambda^*(C(M))\cong \Lambda^*(M) \oplus dt \wedge \Lambda^*(M)$.
Restricting ourselves to the $dt \wedge \Lambda^*(M)$-summand
of this decomposition and applying
 \eqref{_nabla_1_of_rho_Equation_}, we find  
\[
(\nabla_1)_X(t^3\omega\wedge \theta) = t^3 (d\omega \cntrct X)\wedge \theta.
\]
for any $X$ orthogonal to $dt$.
Since $g_1$ is a product metric on $C(M) \cong M\times \R$,
this leads to $\nabla\omega=d\omega$, where
$\nabla$ is the Levi-Civita connection on $M$.
This implies \eqref{_product_conne_to_omega_Equation_}.
We deduced \ref{_NK_from_forms_Theorem_} (iii)
from \ref{_NK_from_forms_Theorem_} (ii).
The proof of \ref{_NK_from_forms_Theorem_} is finished.
\endproof

\hfill

Using \ref{_G_2_unique_Corollary_}, we also obtain the
following useful proposition.

\hfill

\proposition\label{_NK_determined_by_metric_Proposition_}
Let $(M, I, g)$ be a nearly K\"ahler manifold.
Then the almost complex structure is uniquely determined by
the Riemannian structure, unless $M$ is locally isometric to a 
6-sphere.

\endproof


\section[Almost complex structures and their variations] 
{Almost complex structures on 6-manifolds 
and their infinitesimal variations}
\label{_varia_NK_Section_}


\subsection{Hitchin functional and the volume functional}
\label{_extre_volume_via_N_Subsection_}

Let $(M,I)$ be an almost complex 6-manifold with nowhere
degenerate Nijenhuis tensor $N$, and 
$\Vol_I= \det N^* \otimes \overline{\det N^*}$
the corresponding volume form (see
\eqref{_volume_via_N_Equation_}).
In this section we study the extrema of the
functional $I \stackrel \Psi \arrow \int_M \Vol_I$.

A similar functional was studied by N. Hitchin
for 6- and 7-manifolds equipped with a stable 3-form
(see \cite{_Hitchin:stable_}). Since then, this
functional acquired a pivotal role in string theory and 
M-theory, under the name ``Hitchin functional''.

Our first step is to describe the variation of $\Psi$.
We denote by ${\goth M}$ the space of all almost
complex structures with nowhere
degenerate Nijenhuis tensor on $M$.

\hfill

Let $(M, I, \omega)$ be an
almost complex manifold with nowhere degenerate 
Nijenhuis tensor 
\[
N\in \Hom(\Lambda^2 T^{1,0}(M), T^{0,1}(M)),
\]
$\delta\in T_I{\goth M}$ an infinitesimal
variation of $I$, and 
\[ 
   N_\delta\in \Hom(\Lambda^2 T^{1,0}(M), T^{0,1}(M))
\]
the corresponding variation of the Nijenhuis tensor.
Consider the form $\rho:= \omega(N(x,y),z)$
associated with the Hermitian structure
on $M$ as in \ref{_FI_existence_skew_sy_to_Theorem_}.
After rescaling $\omega$, we assume that 
\begin{equation}\label{_rho_constant_Equation_}
   |\rho|_\omega=1.
\end{equation}
Since the Nijenhuis tensor is nowhere degenerate,
$\rho$ is also nowhere degenerate.
Therefore, $\rho$ can be used to identify $T^{0,1}(M)$
and $\Lambda^2 T^{1,0}(M)$, and we may consider $N_\delta$ as 
an endomorphism of $\Lambda^{0,1}(M)$. Notice that this
identification maps $N$ to the identity automorphism of $\Lambda^{0,1}(M)$.

\hfill

\claim \label{_first_vari_Claim_}
In these assumptions,
\begin{equation}\label{_differe_Psi_via_N_Equation_}
\frac{d\Psi}{d I}(\delta) = 2 \Re \int_M \Tr N_\delta\Vol_I.
\end{equation}

\hfill

{\bf Proof:} It is well known that 
\[
\frac{d(\det A)}dt = \det A \Tr\left( A^{-1} \frac{dA}{dt}\right)
\]
for any matrix $A$. Applying that to the map
\[
N^* \otimes \bar N^*:\; 
    \Lambda^{1,0}(M) \otimes \Lambda^{0,1}(M)\arrow
    \Lambda^{0,2}(M) \otimes \Lambda^{2,0}(M),
\]
we obtain that 
\begin{equation}\label{_deri_of_det_of_Vol_I_local_Equation_}
\frac{d(\det(N^* \otimes \bar N^*))}{d I}(\delta)
 = \Tr\left(\frac{(N^*_\delta \otimes \bar N^* + N^* \otimes \bar
  N^*_\delta)} {(N^* \otimes \bar N^*)}\right) 
  \cdot \det(N^*\otimes\bar N^*)
\end{equation}
However, after we identify $\Lambda^{1,0}(M)$
and $\Lambda^{0,2}(M)$ as above, $N$ becomes an identity,
and \eqref{_deri_of_det_of_Vol_I_local_Equation_}
gives 
\begin{equation}\label{_derivative_determ_N_local_Equation_}
\frac{d(\det(N^* \otimes \bar N^*))}{d I}(\delta) = 2  \Re\Tr N_\delta\Vol_I
\end{equation}
\endproof

\hfill

\remark\label{_extre_Psi_via_N_Remark_}
We find that the extrema of the functional
$\Psi(M,I) = \int_M \Vol_I$ are precisely those
almost complex structures for which $\Re \Tr N_\delta=0$
for any infinitesimal variation $\delta$ of $I$.

\subsection{Variations of almost complex structures and
  the Nijenhuis tenor}

It is convenient, following Kodaira and Spencer, to consider
infinitesimal variantions of almost complex structures
as tensors $\delta \in \Lambda^{0,1}(M) \otimes
T^{1,0}(M)$. Indeed, a complex structure on a vector
space $V$, $\dim_\R V =2d$,
can be considered as a point of the Grassmanian
of $d$-dimensional planes in $V\otimes \C$.
The tangent space to a Grassmanian at a point
$W \subset V\otimes \C$ is given by 
$\Hom(W, V\otimes\C/W)$.

Consider the $(0,1)$-part $\nabla^{0,1}$ 
of the Levi-Civita connection
\[
\nabla^{0,1}\delta \in \Lambda^{0,1}(M)  \otimes
T^{1,0}(M)\otimes\Lambda^{0,1}(M),
\]
and let $\bar \6:\; \Lambda^{0,1}(M) \otimes
T^{1,0}(M)\arrow \Lambda^{0,2}(M) \otimes
T^{1,0}(M)$ denote the composiion of
$\nabla^{0,1}$ with the exterior multiplication
map
\[
\Lambda^{0,1}(M)  \otimes
T^{1,0}(M)\otimes\Lambda^{0,1}(M)\arrow \Lambda^{0,2}(M) \otimes
T^{1,0}(M).
\]
The following claim is well known.

\hfill

\claim\label{_Nije_vari_via_conne_Claim_}
Let $(M,I)$ be an almost complex manifold, and 
\[ \delta \in \Lambda^{0,1}(M) \otimes
   T^{1,0}(M)
\] 
an infinitesimal variation
of almost complex structure. Denote by 
$N_\delta\subset \Lambda^{2,0}(M) \otimes
T^{0,1}(M)$ the corresponding
infinitesimal variation of the Nijenhuis tensor
(see Subection \ref{_extre_volume_via_N_Subsection_}).
Then $\bar N_\delta = \bar\6\delta$, where 
\[ \bar \6:\; \Lambda^{0,1}(M) \otimes
   T^{1,0}(M)\arrow \Lambda^{0,2}(M) \otimes
   T^{1,0}(M)
\] is the differential operator defined above.

\hfill

{\bf Proof:} The proof of \ref{_Nije_vari_via_conne_Claim_}
follows from a direct computation (see e.g. 
\cite{_Kodaira_Spencer:Deformation_}).
\endproof

\hfill

\ref{_Nije_vari_via_conne_Claim_}
can be used to study the deformation properties
of the functional $I \stackrel \Psi \arrow \int_M \Vol_I$
constructed above (see Subsection
\ref{_extre_volume_via_N_Subsection_}). Indeed,
from \ref{_extre_Psi_via_N_Remark_} it follows that
$\Psi$ has an extremum at $I$ if and only if 
$\Re\Tr N_\delta=0$ for any $\delta \in \Lambda^{0,1}(M) \otimes
   T^{1,0}(M)$. Using the identification
$T^{1,0}(M)\cong \Lambda^{2,0}(M)$, provided
by the non-degenerate $(3,0)$-form as above,
we can consider $\delta$ as a $(2,1)$-form on $M$.
Then 
\[ \bar\6 \delta\in  
  \Lambda^{0,2}(M)\otimes \Lambda^{2,0}(M)=
  \Lambda^{2,2}(M)
\]
is the  $(2,2)$-part of 
$d \delta$. Under these identifications, and
using $|\rho|_\omega=1$ from \eqref{_rho_constant_Equation_},
we can express
$\Tr \bar N_\delta$ as 
\begin{equation}\label{_Tr_N_delta_Equation_}
\Tr \bar N_\delta = \frac{\bar\6\delta\wedge \omega}{\Vol_I},
\end{equation}
where $\bar\6$ is a $(0,1)$-part of the de Rham differential.
This gives the following claim.

\hfill

\claim\label{_deri_Psi_Claim_}
Let $(M,I, \omega)$ be an almost complex Hermitian
6-manifold with nowhere
degenerate Nijenhuis tensor. Assume that 
 the corresponding 3-form $\rho$ satisfies
$|\rho|_\omega=1$ (see \eqref{_rho_constant_Equation_}). Consider the
functional $\Psi(I)= \int_M \Vol_I$ on the space of such
almost complex structures. Then
\begin{equation}\label{_differe_Psi_via_bar_6_delta_Equation_}
\frac{d\Psi}{d I}(\delta) = 2 \Re\int_M \bar\6\delta\wedge \omega,
\end{equation}
where $\delta\in \Lambda^{0,1}(M) \otimes
   T^{1,0}(M)$ is an infinitesimal deformation of 
an almost complex structure $I$, considered as a $(2,1)$-form on $M$.

\hfill

{\bf Proof:} \ref{_deri_Psi_Claim_} is implied immediately
by \eqref{_Tr_N_delta_Equation_} and \ref{_first_vari_Claim_}.
\endproof

\hfill

Comparing \ref{_deri_Psi_Claim_} with \ref{_extre_Psi_via_N_Remark_},
we find the following

\hfill

\corollary\label{_extremal_I_via_bar_6_Corollary_}
In assumptions of \ref{_deri_Psi_Claim_}, $I$ 
is an extremum of $\Psi$ if and only if 
\begin{equation}\label{_int_bar_6_delta_omega_vanishes_Equation_}
\Re\int_M \bar\6\delta\wedge \omega=0
\end{equation}
for any $\delta\in \Lambda^{2,1}(M)$.

\endproof

\hfill

Integrating by parts, we find that 
\eqref{_int_bar_6_delta_omega_vanishes_Equation_}
is equivalent to 
\[
\Re \int_M \delta\wedge \bar\6\omega=0
\]
and to $\bar\6\omega=0$. This gives the following
proposition

\hfill

\proposition\label{_extre_of_Psi_via_d_omega_Proposition_}
Let $(M,I, \omega)$ be an almost complex 
Hermitian 6-manifold with nowhere
degenerate Nijenhuis tensor.  Consider the
functional $\Psi(I)= \int_M \Vol_I$ on the space of 
such almost complex structures on $M$. Then
$I$ is an extremum of $\Psi$ if and only if
$d\omega$ lies in $\Lambda^{3,0}(M)\oplus \Lambda^{0,3}(M)$.

\endproof

\hfill

Now, \ref{_extre_of_Psi_via_d_omega_Proposition_}
together with \ref{_NK_from_forms_Theorem_}
 implies \ref{_extre_for_M_I_Theorem_}. Notice 
that by \ref{_NK_unique_from_I_Corollary_},
the nearly K\"ahler Hermitian structure on $(M,I)$ 
is (up to a constant multiplier) uniquely determined by $I$.

\hfill

\hfill

{\bf Acknowledgements:} I am grateful to Robert Bryant, Nigel
Hitchin, Paul-Andi Nagy and Uwe Semmelmann for valuable advice and 
consultations. P.-A. Nagy also suggested adding
\ref{_NK_determined_by_metric_Proposition_}.
Much gratitude to the referee for useful suggestions
and the invaluable help of finding a multitude 
of minor errors.

{\small

\hfill

\noindent {\sc Misha Verbitsky\\
University of Glasgow, Department of Mathematics, \\
15 University Gardens, Glasgow G12 8QW, Scotland.}\\
\   \\
{\sc  Institute of Theoretical and
Experimental Physics \\
B. Cheremushkinskaya, 25, Moscow, 117259, Russia }\\
\  \\
\tt verbit@maths.gla.ac.uk, \ \  verbit@mccme.ru 
}

\end{document}